\documentclass[eqthmnum,nocolour,skinny]{jt-calcs}

\usepackage[euler-digits,euler-hat-accent]{eulervm}

\usepackage[backend=bibtex,style=alphabetic,sorting=nyt]{biblatex}
\bibliography{bibliography}

\usepackage{longtable}
\usepackage{pdflscape}
\usepackage{booktabs}

\setlist[enumerate,1]{label={\normalfont(\alph*)}}

\DeclareMathOperator{\Cusp}{Cusp}
\DeclareMathOperator{\Jor}{Jor}
\DeclareMathOperator{\UCh}{UCh}
\newcommand{\cusp}{\mathrm{cusp}}
\newcommand{\rN}{\mathrm{N}}
\newcommand{\rC}{\mathrm{C}}
\newcommand{\rZ}{\mathrm{Z}}
\newcommand{\wc}{\widecheck}
\newcommand{\wt}{\widetilde}

\title{Harish-Chandra Cuspidal Pairs}
\author{Jay Taylor}
\address{Department of Mathematics, The University of Manchester, Oxford Road, Manchester, M13 9PL, United Kingdom.}
\email{jay.taylor@manchester.ac.uk}

\begin{document}
\begin{abstract}
The irreducible characters of a finite reductive group are partitioned into Harish-Chandra series that are labelled by cuspidal pairs. In this note, we describe how one can algorithmically calculate those cuspidal pairs using results of Lusztig.
\end{abstract}

\section{Introduction}
\begin{pa}
The set of (ordinary) irreducible characters $\Irr(\bG^F)$ of a finite reductive group $\bG^F$ is partitioned into Harish-Chandra series
\begin{equation*}
\Irr(\bG^F) = \bigsqcup_{[\bL,\delta] \in \Cusp(\bG^F)/\bG^F} \mathcal{E}(\bG^F, [\bL,\delta]),
\end{equation*}
see \cite[Thm.~5.3.7]{digne-michel:2020:representations-of-finite-groups-of-lie-type} for instance. Here $\Cusp(\bG^F)/\bG^F$ are the orbits, under the natural conjugation action, of \emph{cuspidal pairs} consisting of: a $(\bG,F)$-split Levi subgroup $\bL \leqslant \bG$, meaning $\bL$ is an $F$-stable Levi complement of an $F$-stable parabolic subgroup of $\bG$, and a cuspidal character $\delta \in \Irr(\bL^F)$.

This approach to studying $\Irr(\bG^F)$ has played a crucial role in many recent works solving open problems in the representation theory of finite groups \cite{kessar-malle:2013:quasi-isolated-blocks,malle-spaeth:2016:characters-of-odd-degree,schaeffer-fry:2019:galois-automorphisms,hollenbach:2022:on-e-cuspidal-pairs}. The utility here often comes from the fact that cuspidal characters posses unique features that are not enjoyed by all irreducible characters.

To make this effective one needs to explicitly understand the set $\Cusp(\bG^F)/\bG^F$. In principle this is known by Lusztig's extensive results on the classification of $\Irr(\bG^F)$ \cite{lusztig:1984:characters-of-reductive-groups,lusztig:1988:reductive-groups-with-a-disconnected-centre}. However, in practice, determining exactly which pairs arise requires some work.
\end{pa}

\begin{pa}
As we explain in \cref{sec:jordan-decomp}, see \cref{pa:extraction-process} for a detailed discussion, the work here involves answering the following question: Assume $\bH \leqslant \bG$ is an $F$-stable closed connected subgroup of $\bG$ containing a maximal torus of $\bG$. Then what is the smallest $(\bG,F)$-split Levi subgroup containing $\bH$?

In \cref{sec:split-parabolics} we give an answer to this question in the abstract setting of root systems. Following work of Borel--Tits \cite{borel-tits:1965:groupes-reductifs} and Bonnaf\'e \cite{bonnafe:2006:sln} we give a characterisation of split Levi subsystems in a root system which generalises the usual characterisation found in \cite{bourbaki:2002:lie-groups-chap-4-6}, see \cref{pa:char-para-levi-sys,cor:split-levi-cover}. The connection to the above question is spelled out in \cref{prop:borel-tits}.

In \cite[\S6]{borel-tits:1965:groupes-reductifs} the authors work with a reductive group $\bG$ defined over a field $k$. If $K/k$ is a Galois field extension  then one can define an action of a finite quotient $\Gamma$ of the Galois group $\Gal(K/k)$ on the roots of $\bG$. One can ask a similar question to the one above in this more general setting, and our statements in \cref{sec:split-parabolics} are phrased so that they apply here as well.

Implementing this in Michel's development version of {\sf CHEVIE} \cite{michel:2015:the-development-version-of-CHEVIE} is easy and we have made this available at \cite{taylor:2020:cusppairs}. Thus, we get an algorithmic solution to the above question.
\end{pa}

\begin{pa}
In \cite{broue-malle:1992:theoremes-de-sylow} Brou\'e--Malle introduced the notion of a $d$-split Levi subgroup for an integer $d \geqslant 1$. This generalises the notion of a $(\bG,F)$-split Levi subgroup, with the $(\bG,F)$-split Levi subgroups being the $1$-split Levi subgroups of Brou\'e--Malle. We also generalise our characterisation of split Levi subsystems to this setting, see \cref{cor:char-d-split}. Thus we obtain an algorithmic solution to the above question with ``$(\bG,F)$-split'' replaced by ``$d$-split''.

There is already a lot of functionality within \textsf{CHEVIE} for working with $d$-Harish-Chandra theory. For instance, one can construct all $d$-split Levi subgroups using the command \texttt{SplitLevis} and one can find all $d$-cuspidal unipotent characters using the command \texttt{CuspidalUnipotentCharacters}. Our observations here are concerned with bringing this information for unipotent characters to bear on all cuspidal characters.

The $d$-split Levi subgroups supporting a $d$-cuspidal unipotent character are classified by Brou\'e--Malle--Michel \cite{broue-malle-michel:1993:generic-blocks}. Reducing to the irreducible case, with stored data for the exceptional types computed using the above general \textsf{CHEVIE} functions, we have implemented this classification in \cite{taylor:2020:cusppairs}. With this it is reasonably straightforward to recreate part of the data contained in the tables in \cite{kessar-malle:2013:quasi-isolated-blocks,hollenbach:2022:on-e-cuspidal-pairs}. In these sources one also finds information on the relative Weyl groups, but we shall discuss the computational aspects of this elsewhere.
\end{pa}

\begin{pa}
In \cref{sec:jordan-decomp} we recall results of Lusztig that give a parameterisation of the set of cuspidal pairs $\Cusp(\bG^F)/\bG^F$ in terms of the Jordan decomposition of characters. This reduces us to a question about unipotent characters on the connected centraliser of a semisimple element. As such, we can easily reduce to the setting where our group is adjoint. In \cref{sec:iso-series} we apply the methods of the previous sections when $\bG$ is an adjoint simple group of classical type and the semisimple element is isolated. When $\bG$ is of exceptional type one finds the answer in the tables of \cite{kessar-malle:2013:quasi-isolated-blocks}, or this can be recomputed using \cite{taylor:2020:cusppairs}. We cannot formally reduce to the isolated case but this gives a good flavour of the computations involved.
\end{pa}

\begin{pa}[Notation]
If a group $G$ acts on a set $X$ then we will denote by $X/G$ the resulting set of orbits. Typically the $G$-orbit of $x \in X$ will be denoted by $[x] \in X/G$ and we will sometimes write $G_x$ for the stabiliser of $x$.
\end{pa}

\subsection*{Acknowledgements}
Calculations that culminated in our observations here were originally carried out while the author visited the TU Kaiserslautern. He gratefully acknowledges the support received from grant SFB TRR-195 for this visit as well as the support of the Bell System Fellowship while participating at the IAS. We thank Gunter Malle and Jean Michel for their comments on an earlier version of this work, Ruwen Hollenbach for checking \cite{taylor:2020:cusppairs} against his calculations in \cite{hollenbach:2022:on-e-cuspidal-pairs}, and the referee for their suggestions that improved the exposition of the paper.

\section{Split Parabolic Subgroups of Reflection Groups}\label{sec:split-parabolics}
\begin{pa}\label{pa:pseudo-reflections}
Let $(X,\wc{X})$ be a pair of finite rank free $\mathbb{Z}$-modules equipped with a perfect pairing $\langle -,-\rangle : X \times \wc{X} \to \mathbb{Z}$ and let $\mathbb{K}$ be a field of characteristic $0$. We let $X_{\mathbb{K}} = \mathbb{K} \otimes_{\mathbb{Z}} X$ and $\wc{X}_{\mathbb{K}} = \mathbb{K} \otimes_{\mathbb{Z}} \wc{X}$ be corresponding $\mathbb{K}$-vector spaces and we denote again by $\langle-,-\rangle : X_{\mathbb{K}} \times \wc{X}_{\mathbb{K}} \to \mathbb{K}$ the form obtained from extension by linearity. If $\Sigma \subseteq X_{\mathbb{K}}$ is any subset then we denote by $X_{\mathbb{K},\Sigma} \subseteq X_{\mathbb{K}}$ the $\mathbb{K}$-subspace spanned by $\Sigma$. Recall we have a contravariant bijection $\wc{\phantom{\phi}} : \End_{\mathbb{K}}(X_{\mathbb{K}}) \to \End_{\mathbb{K}}(\wc{X}_{\mathbb{K}})$ such that $\langle \phi x,y\rangle = \langle x,\wc{\phi} y\rangle$ for all $x \in X_{\mathbb{K}}$ and $y \in \wc{X}_{\mathbb{K}}$. For any pair of non-zero vectors $(\alpha,\gamma) \in X \times \wc{X}$ we have an endomorphism $s_{\alpha,\gamma} \in \End_{\mathbb{Z}}(X)$ defined by $s_{\alpha,\gamma}(x) = x-\langle x,\gamma\rangle\alpha$.\end{pa}

\subsection*{Subsets of Root Systems}

\begin{pa}
We assume $\Phi \subseteq X$ is a non-empty finite subset and $\wc{\phantom{\alpha}} : \Phi \to \wc{X}$ is an injection such that $\langle\alpha,\wc{\alpha}\rangle = 2$ and $s_{\alpha,\wc{\alpha}}(\Phi) \subseteq \Phi$ for all $\alpha \in \Phi$. We have $\Phi$ is a root system, in the sense of \cite[VI, \S1.1, Def.~1]{bourbaki:2002:lie-groups-chap-4-6}, in the subspace $X_{\mathbb{R},\Phi}$. We will assume, in addition, that $\Phi$ is reduced. We let $W = \langle s_{\alpha,\wc{\alpha}} \mid \alpha \in \Phi\rangle \leqslant \GL(X_{\mathbb{R}})$ be the corresponding reflection group, which is finite.

Following \cite[VI, \S1.7, Def.~4]{bourbaki:2002:lie-groups-chap-4-6} we say a subset $\Sigma \subseteq \Phi$ is: \emph{closed} if for any $\alpha,\beta \in \Sigma$ with $\alpha+\beta \in \Phi$ we have $\alpha+\beta\in \Sigma$, \emph{symmetric} if $\Sigma = -\Sigma$, and \emph{parabolic} if it is closed and $\Sigma\cup(-\Sigma) = \Phi$. Moreover, we say that $\Sigma \subseteq \Phi$ is a \emph{closed subsystem} if it is closed and symmetric. Such a subsystem $\Sigma \subseteq \Phi$ is naturally a root system in $X_{\mathbb{R},\Sigma}$.
\end{pa}

\begin{pa}\label{pa:char-para-levi-sys}
As $W$ is finite we may choose a $W$-invariant positive definite symmetric bilinear form $(-\mid -) : X_{\mathbb{R}} \times X_{\mathbb{R}} \to \mathbb{R}$. If $\lambda \in X_{\mathbb{R}}$ then the set $\Psi_{\lambda} := \{\alpha \in \Phi \mid (\alpha \mid \lambda) \geqslant 0\}$ is a parabolic set. Moreover, a subset $\Psi \subseteq \Phi$ is parabolic if and only if $\Psi = \Psi_{\lambda}$ for some $\lambda \in X_{\mathbb{R}}$, see \cite[Prop.~3.3.8]{digne-michel:2020:representations-of-finite-groups-of-lie-type}. Given a parabolic set $\Psi \subseteq \Phi$ we have the symmetric part $\Psi_{\mathrm{s}} := \Psi \cap (-\Psi)$ is a closed subsystem, which we call a \emph{Levi subsystem}.

If $E \subseteq X_{\mathbb{R}}$ is a subspace then we have a closed subsystem $\Phi_E := \Phi \cap E \subseteq \Phi$ and by \cite[VI, \S1.7, Prop.~24]{bourbaki:2002:lie-groups-chap-4-6} this is a Levi subsystem. More precisely, a subset $\Sigma \subseteq \Phi$ is a Levi subsystem if and only if $\Sigma = \Phi \cap X_{\mathbb{R},\Sigma}$. Thus, given a subset $\Sigma \subseteq \Phi$ we have $\Phi \cap X_{\mathbb{R},\Sigma}$ is the unique minimal Levi subsystem containing $\Sigma$. We recall the following terminology to be used later.
\end{pa}

\begin{definition}\label{def:Phi-isolated}
We say a subset $\Sigma \subseteq \Phi$ is \emph{$\Phi$-isolated} if $\Phi = \Phi \cap X_{\mathbb{R},\Sigma}$. In other words, $\Sigma$ is not contained in any proper Levi subsystem of $\Phi$.
\end{definition}

\subsection*{Automorphisms}

\begin{pa}\label{pa:perm-action}
Suppose now that $\Gamma \leqslant \rN_{\GL(X_{\mathbb{R}})}(W)$ is finite. As $W\Gamma \leqslant \GL(X_{\mathbb{R}})$ is also finite we may assume our form $(-\mid -)$ is $W\Gamma$-invariant. An element $g \in \Gamma$ must preserve the reflections in $W$, hence permutes the lines spanned by the roots. More precisely, as our root system is reduced we have a well-defined permutation $\rho_g : \Phi \to \Phi$ such that $g\alpha = c_{g,\alpha}\rho_g(\alpha)$ for some real number $c_{g,\alpha} > 0$. We consider $\Phi$ to be a $\Gamma$-set via this permutation action.
\end{pa}

\begin{definition}
A subset $\Sigma \subseteq \Phi$ is said to be a \emph{$\Gamma$-split Levi subsystem} if $\Sigma = \Psi_{\mathrm{s}} = \Psi \cap (-\Psi)$ for some $\Gamma$-stable parabolic subset $\Psi \subseteq \Phi$. If $\phi \in \rN_{\GL(X_{\mathbb{R}})}(W)$ has finite order then we say $\Sigma \subseteq \Phi$ is a $\phi$-split Levi subsystem if it is a $\langle \phi\rangle$-split Levi subsystem.
\end{definition}

\begin{pa}
Our goal now is to prove \cref{cor:split-levi-cover}. If $E \subseteq X_{\mathbb{R}}$ is a $\Gamma$-stable subspace then let $E^{\Gamma}$ be the $\Gamma$-fixed subspace of $E$. If $\theta : X_{\mathbb{R}} \to X_{\mathbb{R}}^{\Gamma}$ is the orthogonal projection onto $X_{\mathbb{R}}^{\Gamma}$ then for $v \in X_{\mathbb{R}}$ we have
\begin{equation}\label{eq:theta}
\theta v = \frac{1}{|\Gamma|}\sum_{g \in \Gamma} gv.
\end{equation}
Hence $(\theta x \mid y) = (\theta x \mid \theta y) = (x\mid\theta y)$ for all $x,y \in X_{\mathbb{R}}$. If $E \subseteq X_{\mathbb{R}}$ is a subset then we let $\Phi^E = \Phi_{E^{\perp}} = \Phi \cap E^{\perp}$, where $E^{\perp} = \{x \in X_{\mathbb{R}} \mid (x\mid e) = 0$ for all $e \in E\}$ is the perpendicular space. If $\lambda \in X_{\mathbb{R}}$ then the symmetric part of the parabolic set $\Psi_{\lambda}$ is $\Phi^{\lambda} := \Phi^{\{\lambda\}}$. As in \cite[Prop.~2.2]{bonnafe:2006:sln} we get the following characterisation of $\Gamma$-split Levi subsystems.
\end{pa}

\begin{prop}\label{prop:split-parabolic-char}
Suppose $\Sigma \subseteq \Phi$ is a subset. Then the following are equivalent:
\begin{enumerate}
	\item $\Sigma$ is a $\Gamma$-split Levi subsystem,
	\item $\Sigma = \Phi^E$ for some subspace $E \subseteq X_{\mathbb{R}}^{\Gamma}$,
	\item $\Sigma = \Phi^{\lambda}$ for some $\lambda \in X_{\mathbb{R}}^{\Gamma}$.
\end{enumerate}
\end{prop}

\begin{proof}
(c) $\Rightarrow$ (b). Clear.

(b) $\Rightarrow$ (a). For $\alpha \in \Phi$ denote by $E^{\alpha} \subseteq E$ the hyperplane $\{e \in E \mid (e\mid\alpha) = 0\}$. By definition $\alpha \in \Phi^E$ if and only if $E^{\alpha} = E$. Now let $F = \bigcup_{\alpha \in \Phi\setminus \Phi^E} E^{\alpha} \subseteq E$, which is empty if $\Phi = \Phi^E$. As $\Phi\setminus \Phi^E$ is finite and $\mathbb{R}$ is infinite it is well known that $E \neq F$ so there exists a vector $\lambda \in E\setminus F$. Taking $\Psi = \Psi_{\lambda}$ we see that $\Phi^E = \Phi^{\lambda} = \Psi_{\mathrm{s}}$ and clearly $\rho_g(\Psi_{\lambda}) = \Psi_{g\lambda} = \Psi_{\lambda}$ for any $g \in \Gamma$.

(a) $\Rightarrow$ (c). As $\Sigma$ is a Levi subsystem we have by \cref{pa:char-para-levi-sys} that $\Sigma = \Phi^{\lambda}$ for some $\lambda \in X_{\mathbb{R}}$. It suffices to show that $\Sigma = \Phi^{\lambda} = \Phi^{\theta\lambda}$. Now, if $\alpha \in \Phi$ then by \cref{eq:theta} we have
\begin{equation*}
(\alpha \mid \theta\lambda) = \frac{1}{|\Gamma|}\sum_{g \in \Gamma} (\alpha \mid g\lambda) = \frac{1}{|\Gamma|}\sum_{g \in \Gamma} c_{g,\alpha}(\rho_g(\alpha) \mid \lambda).
\end{equation*}
If $\alpha \in \Phi^{\lambda}$ then as $\Phi^{\lambda} \subseteq \Psi_{\lambda}$ is $\Gamma$-stable we have $(\rho_g(\alpha) \mid \lambda) \geqslant 0$ for all $g \in \Gamma$ so $(\alpha \mid \theta\lambda) \geqslant 0$, which means $\Phi^{\lambda} \subseteq \Psi_{\theta\lambda}$. As $\Phi^{\lambda}$ is symmetric we must therefore have that $\Phi^{\lambda} \subseteq \Phi^{\theta\lambda}$.

Similarly, if $\alpha \in \Phi\setminus\Phi^{\lambda}$ then as $\Phi\setminus\Phi^{\lambda}$ is $\Gamma$-stable we have $(\rho_g(\alpha) \mid \lambda) < 0$ for all $g \in \Gamma$ so $(\alpha \mid \theta\lambda) < 0$. Hence, $\Phi\setminus\Phi^{\lambda} \subseteq \Phi\setminus \Psi_{\theta\lambda}$ and as $\Phi\setminus\Phi^{\lambda}$ is symmetric we must have $\Phi\setminus\Phi^{\lambda} \subseteq \Phi\setminus\Phi^{\theta\lambda}$.
\end{proof}

\begin{lem}\label{lem:split-levi}
If $E \subseteq X_{\mathbb{R}}^{\Gamma}$ is a subspace then $\Phi_{\theta^{-1}(E)}$ is a $\Gamma$-split Levi subsystem of $\Phi$.
\end{lem}

\begin{proof}
By the equivalence of (a) and (b) in \cref{prop:split-parabolic-char} it suffices to show that $\theta^{-1}(E) = \theta(E^{\perp})^{\perp}$ as then $\Phi_{\theta^{-1}(E)} = \Phi^{\theta(E^{\perp})}$. We have $\theta^{-1}(E) \subseteq \theta(E^{\perp})^{\perp}$ because if $v \in \theta^{-1}(E)$ and $x \in E^{\perp}$ then $(v\mid \theta x) = (\theta v \mid x) = 0$. Conversely, by assumption we have $(X_{\mathbb{R}}^{\Gamma})^{\perp} \subseteq E^{\perp}$ which implies that $\theta(E^{\perp}) = E^{\perp} \cap X_{\mathbb{R}}^{\Gamma}$ and so $\theta(E^{\perp})^{\perp} = (E^{\perp})^{\perp} + (X_{\mathbb{R}}^{\Gamma})^{\perp} = E \oplus (X_{\mathbb{R}}^{\Gamma})^{\perp}$. Hence $\theta(\theta(E^{\perp})^{\perp}) = E$ which gives $\theta^{-1}(E) = \theta(E^{\perp})^{\perp}$.
\end{proof}

\begin{cor}\label{cor:split-levi-cover}
For any subset $\Sigma \subseteq \Phi$ we have $\Phi \cap \theta^{-1}(\theta(X_{\mathbb{R},\Sigma}))$ is the smallest $\Gamma$-split Levi subsystem containing $\Sigma$. In particular, $\Sigma$ is a $\Gamma$-split Levi subsystem if and only if $\Sigma = \Phi \cap \theta^{-1}(\theta(X_{\mathbb{R},\Sigma}))$.
\end{cor}

\begin{proof}
By \cref{lem:split-levi} we certainly have $\Phi \cap \theta^{-1}(\theta(X_{\mathbb{R},\Sigma}))$ is a $\Gamma$-split Levi subsystem containing $\Sigma$. However, if $\Psi \supseteq \Sigma$ is a $\Gamma$-split Levi subsystem then $X_{\mathbb{R},\Sigma} \subseteq X_{\mathbb{R},\Psi}$ which implies that $\theta^{-1}(\theta(X_{\mathbb{R},\Sigma})) \subseteq \theta^{-1}(\theta(X_{\mathbb{R},\Psi}))$. So $\Phi \cap \theta^{-1}(\theta(X_{\mathbb{R},\Sigma})) \subseteq \Psi = \Phi \cap \theta^{-1}(\theta(X_{\mathbb{R},\Psi}))$.
\end{proof}

\begin{pa}\label{pa:code-example}
In \textsf{CHEVIE} \cite{michel:2015:the-development-version-of-CHEVIE} one finds the command \texttt{EigenspaceProjector} which can be used to calculate the projector onto the $\phi$-fixed subspace of $X_{\mathbb{R}}$ for any finite order element $\phi \in \rN_{\GL(X_{\mathbb{R}})}(W)$. With this it is easy to find those roots $\alpha \in \Phi$ satisfying $\theta\alpha \in \theta(X_{\mathbb{R},\Sigma}) = X_{\mathbb{R},\theta(\Sigma)}$, see the function \texttt{SplitLeviCover} in \cite{taylor:2020:cusppairs} for more details.
\end{pa}

\subsection{Positive Roots}
\begin{pa}\label{pa:standard-parabolic}
Given $\phi \in \rN_{\GL(X_{\mathbb{R}})}(W)$ let $\mathcal{L}(\Phi,\phi)$ denote the set of pairs $(\Sigma,w)$ with $\Sigma \subseteq \Phi$ a $w\phi$-split Levi subsystem. There is a natural action of $W$ on $\mathcal{L}(\Phi,\phi)$ given by
\begin{equation*}
x \cdot (\Sigma,w) = (x\Sigma,xw \iota_{\phi}(x)^{-1}),
\end{equation*}
where $\iota_{\phi} : W \to W$ is the automorphism defined by $\iota_{\phi}(y) = \phi y \phi^{-1}$.

Assume $\Phi^+ \subseteq \Phi$ is a fixed system of positive roots. If $g \in \rN_{\GL(X_{\mathbb{R}})}(W)$ then $\rho_g(\Phi^+) \subseteq \Phi$ is another positive system and as $W$ acts simply transitively on all such positive systems we see that there exists a unique element $\phi \in Wg$ satisfying $\rho_{\phi}(\Phi^+) = \Phi^+$. Note that in the discussion leading up to \cref{cor:split-levi-cover} we did not make reference to any such choice of positive roots.

Now let $\Delta \subseteq \Phi^+$ be the unique simple system of roots contained in $\Phi^+$. Each subset $I \subseteq \Delta$ gives rise to a Levi subsystem $\Phi_I := \Phi_{X_{\mathbb{R},I}} = \Phi \cap X_{\mathbb{R},I}$ which is the symmetric part of the parabolic set $\Psi_I = \Phi_I\cup \Phi^+$. Both $\Phi_I$ and $\Psi_I$ are often called \emph{standard}. The corresponding reflection group $W_I = \langle s_{\alpha,\wc\alpha} \mid \alpha \in I\rangle \leqslant W$ is a standard parabolic subgroup. The following shall be utilised in the next section.
\end{pa}

\begin{prop}\label{lem:conj-to-standard}
Assume $\phi \in \rN_{\GL(X_{\mathbb{R}})}(W)$ satisfies $\rho_{\phi}(\Phi^+) = \Phi^+$. Then each orbit $\mathcal{L}(\Phi,\phi)/W$ contains a pair of the form $(\Phi_I,z)$, where $\Phi_I$ is a $\phi$-split standard Levi subsystem with $I \subseteq \Delta$ and $z \in W_I$.
\end{prop}

\begin{proof}
Assume $(\Sigma,w) \in \mathcal{L}(\Phi,\phi)$ and let $\Gamma = \langle w\phi\rangle$. If $V = X_{\mathbb{R},\Phi}$ then we have a $\Gamma$-stable orthogonal decomposition $X_{\mathbb{R}} = V\oplus V^{\perp}$. It thus follows easily from \cref{prop:split-parabolic-char} that there exists a vector $\lambda \in V^{\Gamma}$ such that $\Sigma = \Phi^{\lambda}$.

The set $D = \{x \in V \mid (x\mid \alpha) \geqslant 0$ for all $\alpha \in \Delta\}$ is a fundamental domain for the action of $W$, see \cite[VI, \S1.5, Thm.~2]{bourbaki:2002:lie-groups-chap-4-6}. Hence, there exists an element $x \in W$ such that $\mu = x\lambda \in D$. If $z = xw\iota_{\phi}(x)^{-1}$ then $\mu$ is fixed by $z\phi$. Now $\phi D = D$, because $\rho_{\phi}(\Delta) = \Delta$, so $\phi\mu \in D$ which means $\phi\mu = \mu$ and $z\mu = \mu$ because $z\phi\mu = \mu$. We then have $x\Sigma = \Phi^{\mu} = \Phi_I$, where $I = \{\alpha \in \Delta \mid (\alpha \mid \mu) = 0\}$, and $z \in W_I$, see \cite[Cor.~A.28]{malle-testerman:2011:linear-algebraic-groups} for instance.
\end{proof}

\subsection{A Generalisation}
\begin{pa}
We now consider a generalisation of these ideas following Brou\'e--Malle \cite{broue-malle:1992:theoremes-de-sylow}. For this we assume $\phi \in \rN_{\GL(X)}(W)$ is an automorphism of finite order $n > 0$, which we consider as an element of $\GL(X_{\mathbb{K}})$ in the natural way. In other words, we assume that $\phi$ is defined over $\mathbb{Z}$. For an integer $d > 0$ we denote by $\mu_d \leqslant \mathbb{C}^{\times}$ the subgroup of $d$th roots of unity and by $\mu_d^* \subseteq \mu_d$ the subset of elements of order $d$, i.e., the primitive $d$th roots of unity.

Assume $\mathbb{K \subseteq C}$ is a subfield and $E \subseteq X_{\mathbb{K}}$ is a $\phi$-stable $\mathbb{K}$-subspace. If $\zeta \in \mathbb{K}^{\times}$ then we denote by $E(\phi,\zeta) \subseteq E$ the $\zeta$-eigenspace of $\phi$ acting on $E$. Moreover, for an integer $d > 0$ we let
\begin{equation*}
E(\phi,d) = \bigoplus_{\zeta \in \mu_d^*} E(\phi,\zeta).
\end{equation*}
Of course, $E(\phi,d) = \{0\}$ if $d \nmid n$ and $E = \bigoplus_{d \mid n} E(\phi,d)$ if $\mu_n \subseteq \mathbb{K}$. From now on $\mathbb{K} = \mathbb{Q}(\mu_n) \subseteq \mathbb{C}$ is the $n$th cyclotomic field and $\mathcal{G} = \Gal(\mathbb{K/Q})$.

We have a natural semilinear action of $\mathcal{G}$ on $X_{\mathbb{K}} = \mathbb{K} \otimes_{\mathbb{Q}} X_{\mathbb{Q}}$, in the sense of \cite[V, \S10.4]{bourbaki:2003:algebra-II-4-7}, where $\sigma \in \mathcal{G}$ acts via $\sigma\otimes_{\mathbb{Q}}\Id$ on $X_{\mathbb{K}}$. This action commutes with the action of $\phi$ so it is clear that $\sigma X_{\mathbb{K}}(\phi,\zeta) = X_{\mathbb{K}}(\phi,\sigma(\zeta))$ for any $\sigma \in \mathcal{G}$. Hence the space $X_{\mathbb{K}}(\phi,d)$ is $\mathcal{G}$-stable so there is a canonical $\mathbb{Q}$-subspace $X_{\mathbb{Q}}(\phi,d) := X_{\mathbb{K}}(\phi,d)^{\mathcal{G}}$, the subspace of $\mathcal{G}$-invariants, satisfying $X_{\mathbb{K}}(\phi,d) = \mathbb{K} \otimes_{\mathbb{Q}} X_{\mathbb{Q}}(\phi,d)$, see \cite[V, \S10.4, Cor.]{bourbaki:2003:algebra-II-4-7}.
\end{pa}

\begin{pa}
Fix an indeterminate $\mathbf{q}$. For an integer $d > 0$ we denote by $\Phi_d(\mathbf{q})$ the $d$th cyclotomic polynomial $\prod_{\zeta \in \mu_d^*} {(\mathbf{q}-\zeta)} \in \mathbb{Z}[\mathbf{q}]$. We have $X_{\mathbb{Q}} = \bigoplus_{d \mid n} X_{\mathbb{Q}}(\phi,d)$ and the characteristic polynomial of $\phi$ on $X_{\mathbb{Q}}(\phi,d)$ is given by $\Phi_d(\mathbf{q})^{a_d}$, where $a_d = \dim_{\mathbb{K}}(X_{\mathbb{K}}(\phi,\zeta))$ with $\zeta \in \mu_d^*$. Here we use that $\mathcal{G}$ acts transitively on $\mu_d^*$, see \cite[V, \S11.5, Thm.~2]{bourbaki:2003:algebra-II-4-7}. If $E \subseteq X_{\mathbb{Q}}(\phi,d)$ is a $\phi$-stable $\mathbb{Q}$-subspace then the characteristic polynomial of $\phi$ on $E$ is given by $\Phi_d(\mathbf{q})^b$ for some $0 \leqslant b \leqslant a_d$.

We assume $(-\mid -)$ is a non-degenerate $W\langle\phi\rangle$-invariant positive definite symmetric bilinear form $X_{\mathbb{Q}}$. We then consider $(-\mid-)$ as a bilinear form on $X_{\mathbb{K}}$ by scalar extension. As before, $E^{\perp} \subseteq X_{\mathbb{K}}$ is the perpendicular space of $E \subseteq X_{\mathbb{K}}$ and $\Phi^E = \Phi \cap E^{\perp}$. After \cite{broue-malle:1992:theoremes-de-sylow} we make the following definition.
\end{pa}

\begin{definition}
If $d > 0$ is an integer then a subset $\Sigma \subseteq \Phi$ is a \emph{$(\phi,d)$-split Levi subsystem} if $\Sigma = \Phi^E$ for some $\phi$-stable $\mathbb{Q}$-subspace $E \subseteq X_{\mathbb{Q}}(\phi,d)$.
\end{definition}

\begin{lem}\label{lem:def-cmplx-split-levi}
Recall our assumption that $\phi \in \rN_{\GL(X)}(W)$ is defined over $\mathbb{Z}$. Fix a root of unity $\zeta \in \mu_d^*$. For any subset $\Sigma \subseteq \Phi$ the following are equivalent:
\begin{enumerate}
	\item $\Sigma$ is $(\phi,d)$-split,
	\item $\Sigma = \Phi^E = \Phi \cap E^{\perp}$ for some subspace $E \subseteq X_{\mathbb{K}}(\phi,\zeta)$, which is necessarily $\phi$-stable.
\end{enumerate}
\end{lem}

\begin{proof}
Let $E \subseteq X_{\mathbb{Q}}(\phi,d)$ be a $\phi$-stable $\mathbb{Q}$-subspace and set $E_{\mathbb{K}} = \mathbb{K} \otimes_{\mathbb{Q}} E$ and $\tilde{E} = E_{\mathbb{K}}(\phi,\zeta)$. It suffices to show that $\Phi^E = \Phi^{\tilde{E}}$. As above, we have $\mathcal{G}$ acts on $E_{\mathbb{K}}$ and $E_{\mathbb{K}} = \bigoplus_{\sigma \in \mathcal{G}/\mathcal{G}_{\zeta}} \sigma \tilde{E}$. However, if $\alpha \in \Phi$ and $e \in \tilde{E}$ then $\sigma (\alpha \mid e) = (\alpha \mid \sigma e)$. This shows that $\Phi^{\tilde{E}} = \Phi^{E_{\mathbb{K}}} = \Phi^E$.
\end{proof}

\begin{pa}
For $\zeta \in \mathbb{K}^{\times}$ we denote by $\theta_{\zeta} : X_{\mathbb{K}} \to X_{\mathbb{K}}(\phi,\zeta)$ the orthogonal projection onto the $\zeta$-eigenspace of $\phi$. The exact same proofs used to obtain \cref{lem:split-levi,cor:split-levi-cover} yield the following, which give a characterisation of $(\phi,d)$-split Levi subsystems.
\end{pa}

\begin{lem}
If $\zeta \in \mu_d^*$ is a primitive $d$th root of unity and $E \subseteq X_{\mathbb{K}}(\phi,\zeta)$ is a $\mathbb{K}$-subspace then $\Phi \cap \theta_{\zeta}^{-1}(E)$ is a $(\phi,d)$-split Levi subsystem of $\Phi$.
\end{lem}

\begin{cor}\label{cor:char-d-split}
Let $\zeta \in \mu_d^*$ be a primitive $d$th root of unity. Then for any subset $\Sigma \subseteq \Phi$ we have $\Phi \cap \theta_{\zeta}^{-1}(\theta_{\zeta}(X_{\mathbb{K},\Sigma}))$ is the smallest $(\phi,d)$-split Levi subsystem containing $\Sigma$. In particular, $\Sigma$ is a $(\phi,d)$-split Levi subsystem if and only if $\Sigma = \Phi \cap \theta_{\zeta}^{-1}(\theta_{\zeta}(X_{\mathbb{K},\Sigma}))$.
\end{cor}

\begin{pa}\label{pa:exmp-d-split}
The eigenspace projector $\theta_{\zeta}$ can also be obtained through the \textsf{CHEVIE} command \texttt{EigenspaceProjector} and so the function \texttt{SplitLeviCover} provided in \cite{taylor:2020:cusppairs} will also deal with this more general setting. Now suppose $\phi \in \rN_{\GL(X_{\mathbb{R}})}(W)$ has finite order $n > 0$, so $\phi$ does not necessarily stabilise $X$. Jean Michel has informed us that in this situation one simply takes (b), with $\mathbb{K} = \mathbb{Q}(\mu_n)$, of \cref{lem:def-cmplx-split-levi} as the definition for $(\phi,d)$-split. With this definition \cref{cor:char-d-split} also applies in this situation.
\end{pa}

\section{Jordan Decomposition and Cuspidal Pairs}\label{sec:jordan-decomp}
\begin{pa}\label{pa:reductive-setup}
From now on $\bG$ is a connected reductive algebraic group defined over an algebraic closure $\mathbb{F} = \overline{\mathbb{F}_p}$ of the finite field of prime order $p>0$. We let $F : \bG \to \bG$ be a Frobenius root as in \cite{digne-michel:2020:representations-of-finite-groups-of-lie-type}, by which we mean $F^n$ is a Frobenius endomorphism for some integer $n \geqslant 1$. The set of unipotent characters of the finite group $\bG^F$ will be denoted by $\UCh(\bG,F)$.

We fix an $F$-stable maximal torus and Borel subgroup $\bT_0 \leqslant \bB_0$ of $\bG$. The roots of $\bG$ are denoted by $\Phi \subseteq X = X(\bT_0) = \Hom(\bT_0,\mathbb{F}^{\times})$ and the Weyl group by $W = \rN_{\bG}(\bT_0)/\bT_0$. Taking $\wc{X} = \Hom(\mathbb{F}^{\times},\bT_0)$ we have a natural perfect pairing $\langle -,-\rangle : X \times \wc{X} \to \mathbb{Z}$ that puts us in the setting of \cref{sec:split-parabolics}. Recall that $\bB_0$ determines a set of simple and positive roots $\Delta \subseteq \Phi^+ \subseteq \Phi$. We let $\overline{\phantom{x}} : \rN_{\bG}(\bT_0) \to W$ be the natural quotient map.

The endomorphism of $X$ given by $\chi \mapsto \chi\circ F$ will again be denoted by $F$. Note this assignment is contravariant. Extending linearly this defines an endomorphism of $X_{\mathbb{R}}$ that factorises uniquely as $q\phi$ with $q \in \mathbb{R}_{>0}$ a positive real number and $\phi \in \GL(X_{\mathbb{R}})$ of finite order. The automorphism $\phi$ satisfies $\rho_{\phi}(\Phi^+) = \Phi^+$ as in \cref{pa:standard-parabolic}.

For an element $g \in \bG$ we let $\iota_g : \bG \to \bG$ denote the inner automorphism given by $\iota_g(x) = gxg^{-1}$ and we let $gF := \iota_g\circ F$.
\end{pa}

\begin{pa}
Consider the subset $\bG_{\sss} \subseteq \bG$ of semisimple elements. The set $\mathfrak{T}(\bG,F) = \{(s,n) \in \bG_{\sss} \times \bG \mid {}^nF(s) = s\}$ is a $\bG$-set with the action given by
\begin{equation}\label{eq:G-action}
g\cdot (s,n) = ({}^gs, gnF(g)^{-1}).
\end{equation}
Note the stabiliser of $(s,n)$ under this action is the finite group $\rC_{\bG}(s)^{nF}$. By the Lang--Steinberg Theorem we have a bijective map $\bG_{\sss}^F/\bG^F \to \mathfrak{T}(\bG,F)/\bG$ given by $[s] \mapsto [s,1]$, where $\bG_{\sss}^F/\bG^F$ is the set of semisimple conjugacy classes of $\bG^F$.

For $s \in \bT_0$ we let $\rC_W(s) := \rN_{\rC_{\bG}(s)}(\bT_0)/\bT_0$ and $\rC_W^{\circ}(s) = \rN_{\rC_{\bG}^{\circ}(s)}(\bT_0)/\bT_0$, where $\rC_{\bG}^{\circ}(s)$ is the connected component of the centraliser $\rC_{\bG}(s)$. Let $\mathcal{A}_W(s,F)$ be the set of cosets $a \in \rC_W^{\circ}(s) \backslash W$ such that ${}^wF(s) = s$ for some (any) $w \in a$. Then $\mathcal{T}_{\bG}(\bT_0,F)$ denotes the set of pairs $(s,a)$ with $s \in \bT_0$ and $a \in \mathcal{A}_W(s,F)$. There is a natural action of $W$ on $\mathcal{T}_{\bG}(\bT_0,F)$, defined exactly as in \cref{eq:G-action}. We remark that we have the following well-known parameterisation of the semisimple classes of $\bG^F$.
\end{pa}

\begin{lem}\label{lem:ss-reps}
We have a well-defined bijection $\mathcal{T}_{\bG}(\bT_0,F)/W \to \mathfrak{T}(\bG,F)/\bG$ defined by sending $[s,\rC_W^{\circ}(s)\overline{n}] \mapsto [s,n]$, where $n \in \rN_{\bG}(\bT_0)$.
\end{lem}

\begin{pa}\label{pa:induced-Frob-map}
Let $\bS \leqslant \bG$ be an $F$-stable torus. We say $\bS$ is $F$-split if the endomorphism of $X(\bS)_{\mathbb{R}} = \mathbb{R}\otimes_{\mathbb{Z}} X(\bS)$ induced by $F$ is multiplication by $q \in \mathbb{R}_{>0}$. Any $F$-stable subtorus of $\bG$ contains a unique maximal $F$-split subtorus. We denote by $\rZ_F^{\circ}(\bG)$ this subtorus of $\rZ^{\circ}(\bG)$. The following result is well known. We give a few details just to make the connection with the material in \cref{sec:split-parabolics}.
\end{pa}

\begin{prop}[Borel--Tits]\label{prop:borel-tits}
Suppose $\bH \leqslant \bG$ is an $F$-stable closed connected maximal rank subgroup of $\bG$. Then there is a unique minimal $(\bG,F)$-split Levi subgroup containing $\bH$ and this subgroup is $\rC_{\bG}(\rZ_F^{\circ}(\bH))$. In particular, $\bH$ is a $(\bG,F)$-split Levi subgroup if and only if $\bH = \rC_{\bG}(\rZ_F^{\circ}(\bH))$.
\end{prop}

\begin{proof}
By assumption we have ${}^g\bT_0 \leqslant \bH$ for some $g \in \bG$ with $g^{-1}F(g) = n \in \rN_{\bG}(\bT_0)$. Hence, we can assume that $\bH \geqslant \bT_0$ if we replace $F$ by $nF$. Let $q\tau$ be the factorisation of $nF$ on $X_{\mathbb{R}}$ with $\tau \in \GL(X_{\mathbb{R}})$ of finite order. As in \cite[\S3.3]{digne-michel:2020:representations-of-finite-groups-of-lie-type} we have a map $\Sigma \mapsto \bG_{\Sigma}$ which gives a bijection between the lattice of quasi-closed subsets of $\Phi$ and the closed connected subgroups of $\bG$ containing $\bT_0$. If $\Gamma = \langle \tau \rangle$ then this maps the $\Gamma$-stable parabolic subsets of $\Phi$ onto the $F$-stable parabolic subgroups of $\bG$ containing $\bT_0$, see \cite[Prop.~3.4.5]{digne-michel:2020:representations-of-finite-groups-of-lie-type}, and thus the $\Gamma$-split Levi subsystems of $\Phi$ onto the $F$-stable Levi subgroups of $\bG$ containing $\bT_0$. The first statement follows from \cref{cor:split-levi-cover}, the second from \cite[Lem.~3.1.3]{geck-malle:2020:the-character-theory-of-finite-groups-of-lie-type}, and the last statement is obvious.
\end{proof}

\subsection*{Jordan Parameters}
\begin{pa}\label{pa:action-jordan}
We use a language for Jordan decompositions similar to that used by Cabanes--Sp\"ath \cite[\S8]{cabanes-spaeth:2015:equivariance-and-extendibility-in-groups-of-type-A}. If $(s_1,n_1),(s_2,n_2) \in \mathfrak{T}(\bG,F)$ are such that $(s_2,n_2) = ({}^gs_1,gn_1F(g)^{-1})$, for some $g \in \bG$, then we have $\iota_g\circ n_1F = n_2F \circ \iota_g$. We define $\Jor^{\circ}(\bG,F)$ to be the set of pairs $((s,n),\psi)$ with $(s,n) \in \mathfrak{T}(\bG,F)$ and $\psi \in \UCh(\rC_{\bG}^{\circ}(s),nF)$. In other words, we have
\begin{equation*}
\Jor^{\circ}(\bG,F) = \bigsqcup_{(s,n) \in \mathfrak{T}(\bG,F)} \UCh(\rC_{\bG}^{\circ}(s),nF).
\end{equation*}
Moreover, $\bG$ acts on $\Jor^{\circ}(\bG,F)$ via $g\cdot (\mathfrak{s},\psi) = (g\cdot \mathfrak{s},{}^g\psi)$ where ${}^g\psi = \psi\circ\iota_g^{-1}$.

Now let $\bG_{\ad} = \bG/\rZ(\bG)$, which for these purposes we consider to be the adjoint group of $\bG$. We denote again by $F$ the Frobenius root of $\bG_{\ad}$ induced by that of $\bG$. If $g\rZ(\bG) \in \bG_{\ad}^F$ then for any $x \in \bG^F$ we have ${}^gx \in \bG^F$ and this gives a well-defined action of $\bG_{\ad}^F$ on $\bG^F$, hence also of $\bG_{\ad}^F$ on $\Irr(\bG^F)$.

Recall from the introduction that $\Cusp(\bG^F)$ is the set of pairs $(\bL,\delta)$ with $\bL \leqslant \bG$ a $(\bG,F)$-split Levi subgroup and $\delta \in \Irr(\bL^F)$ a cuspidal character. We get an action of $\bG_{\ad}^F$ on $\Cusp(\bG,F)$ by setting $g\rZ(\bG) \cdot (\bL,\delta) = ({}^g\bL,{}^g\delta)$ where ${}^g\delta = \delta\circ\iota_g^{-1}$.
\end{pa}

\begin{pa}\label{pa:compatible-Jordan}
We now fix a tuple $(\bG^{\star},\bB_0^{\star},\bT_0^{\star},F^{\star})$ that is dual to $(\bG,\bB_0,\bT_0,F)$. We will denote by $W^{\star} = \rN_{\bG^{\star}}(\bT_0^{\star})/\bT_0^{\star}$ the dual Weyl group and by $\Delta^{\star} \subseteq \Phi^{\star+} \subseteq \Phi^{\star} \subseteq X(\bT_0^{\star})$ the corresponding sets of roots.

Recall that in \cite{lusztig:1984:characters-of-reductive-groups,lusztig:1988:reductive-groups-with-a-disconnected-centre} Lusztig has shown the existence of a \emph{Jordan decomposition} $J : \Irr(\bG^F) \to \Jor^{\circ}(\bG^{\star},F^{\star})/\bG^{\star}$ which is a certain surjective map whose fibres are the orbits $\Irr(\bG^F)/\bG_{\ad}^F$. In what follows we will need to assume some compatibility properties between Jordan decompositions. We will write $J^{\bullet}$ to denote a family of Jordan decompositions $J^{\bL} : \Irr(\bL^F) \to \Jor^{\circ}(\bL^{\star},F^{\star})/\bL^{\star}$, one for each $F$-stable Levi subgroup $\bL \leqslant \bG$. Implicitly this involves the choice of an $F^{\star}$-stable Levi subgroup $\bL^{\star} \leqslant \bG^{\star}$ corresponding to $\bL$ under duality, see \cite[\S11.4]{digne-michel:2020:representations-of-finite-groups-of-lie-type}.

We say the family $J^{\bullet}$ is \emph{$\bG^F$-invariant} if for each $F$-stable Levi subgroup $\bL \leqslant \bG$ and element $g \in \bG^F$ the following hold:
\begin{itemize}
	\item $({}^g\bL)^{\star} = \bL^{\star}$
	\item $J^{\bL}(\chi) = J^{{}^g\bL}({}^g\chi)$ for all $\chi \in \Irr(\bL^F)$.
\end{itemize}
Requiring these properties is equivalent to choosing one map $J^{\bL}$ for each $\bG^F$-orbit of Levi subgroups.
\end{pa}

\begin{rem}
Note this assumption implies that if $\bL_1,\bL_2 \leqslant \bG$ are two $F$-stable Levi subgroups such that $\bL_1^{\star} = {}^g\bL_2^{\star}$ for some $g \in \bG^{\star F^{\star}}$ then $\bL_1^{\star} = \bL_2^{\star}$ and $\bL_1$ is $\bG^F$-conjugate to $\bL_2$.
\end{rem}

\subsection{Cuspidal Characters}
\begin{pa}
We start by reducing our problem concerning cuspidal pairs to a question in the dual group using the Jordan decomposition. For this, let us define the set $\mathcal{J}_{\cusp}^{\circ}(\bG,F)$ of tuples $((s,n),(\bL,\psi))$ such that:
\begin{itemize}
	\item $(s,n) \in \mathfrak{T}(\bG,F)$ with $s \in \bL$,
	\item $\bL \leqslant \bG$ is the smallest $(\bG,nF)$-split Levi subgroup containing $\rC_{\bL}^{\circ}(s)$,
	\item $\psi \in \UCh(\rC_{\bL}^{\circ}(s),nF)$ is a cuspidal character.
\end{itemize}

Note that as $s \in \bL$ we have $\rC_{\bL}^{\circ}(s)$ contains a maximal torus of $\bL$, hence $\bG$, so \cref{prop:borel-tits} applies. In particular, the second condition can equivalently be stated as $\rZ_{nF}^{\circ}(\bL) = \rZ_{nF}^{\circ}(\rC_{\bL}^{\circ}(s))$. Of course, there is an action of $\bG$ on $\mathcal{J}_{\cusp}^{\circ}(\bG,F)$ given by $g\cdot (\mathfrak{s},(\bL,\psi)) = (g\cdot\mathfrak{s},({}^g\bL,{}^g\psi))$. With this we have the following consequence of Lusztig's characterisation of cuspidal characters given in \cite[2.18]{lusztig:1978:representations-of-finite-chevalley-groups} and \cite[7.8]{lusztig:1977:irreducible-representations-of-finite-classical-groups}.
\end{pa}

\begin{thm}[Lusztig]\label{thm:jor-dec-cusp}
Assume $J^{\bullet}$ is a $\bG^F$-invariant family of Jordan decompositions as in \cref{pa:compatible-Jordan}. Then we have a well-defined bijection
\begin{equation}
\begin{aligned}
\Cusp(\bG^F)/\bG_{\ad}^F &\to \mathcal{J}_\cusp^{\circ}(\bG^{\star},F^{\star})/\bG^{\star}\\
[\bL,\delta] &\mapsto [\mathfrak{s},(\bL^{\star},\psi)],
\end{aligned}
\end{equation}
where $J^\bL(\delta) = [\mathfrak{s},\psi]$.
\end{thm}

\begin{proof}
Assume $(\bL,\delta) \in \Cusp(\bG^F)$ and $J^{\bL}(\delta) = [\mathfrak{s},\psi]$ then by \cite[Thm.~3.2.22]{geck-malle:2020:the-character-theory-of-finite-groups-of-lie-type} we have $(\mathfrak{s},(\bL^{\star},\psi)) \in \mathcal{J}_{\cusp}^{\circ}(\bG^{\star},F^{\star})$. In particular, this map makes sense. To see it is surjective simply note that we can choose a tuple of the form $((s,1),(\bL^{\star},\psi)) \in \mathcal{J}_{\cusp}^{\circ}(\bG^{\star},F^{\star})$ in each orbit. There then exists a character $\delta \in \Irr(\bL^F)$ satisfying $J^{\bL}(\delta) = [(s,1),\psi]$ and any such $\delta$ is cuspidal, again by \cite[Thm.~3.2.22]{geck-malle:2020:the-character-theory-of-finite-groups-of-lie-type}.

Let $(\bL_1,\delta_1), (\bL_2,\delta_2) \in \Cusp(\bG^F)$. Assume first that $(\bL_2,\delta_2) = ({}^g\bL_1,{}^g\delta_1)$ for some $g\rZ(\bG) \in \bG_{\ad}^F$. As $\rZ(\bG) \leqslant \bL_1$ we have by the Lang--Steinberg Theorem that $g^{-1}F(g) = l^{-1}F(l)$ for some $l \in \bL_1$ so $\bL_2 = {}^h\bL_1$ with $h := gl^{-1} \in \bG^F$. Moreover, $l\rZ(\bL) \in \bL_{\ad}^F$ because $\rZ(\bG) \leqslant \rZ(\bL)$ so $J^{\bL_2}(\delta_2) = J^{{}^h\bL_1}({}^{hl}\delta_1) = J^{\bL_1}({}^l\delta_1) = J^{\bL_1}(\delta_1)$ by our assumption that $J^{\bullet}$ is $\bG^F$-invariant. Thus, the map is well defined.

Now let $J^{\bL_1}(\delta_1) = [(s_1,1),\psi_1]$ and $J^{\bL_2}(\delta_2) = [(s_2,1),\psi_2]$ and assume that we have $((s_1,1),(\bL_1^{\star},\psi_1)) = x\cdot ((s_2,1),(\bL_2^{\star},\psi_2))$ for some $x \in \bG^{\star}$. Note $x \in \bG^{\star F^{\star}}$ so $\bL_1^{\star} = \bL_2^{\star}$ and $\bL_1 = {}^g\bL_2$ for some $g \in \bG^F$. Now we have $J^{\bL_1}(\delta_1) = [(s_1,1),\psi_1] = [(s_2,1),\psi_2] = J^{\bL_2}(\delta_2) = J^{\bL_1}({}^g\delta_2)$. It follows that ${}^{lg}\delta_2 = \delta_1$ for some element $l\rZ(\bL_1) \in (\bL_1/\rZ(\bL_1))^F$.

By \cite[Prop.~4.2]{bonnafe:2006:sln} we have $\rZ(\bL_1) = \rZ^{\circ}(\bL_1)\rZ(\bG)$ so by the Lang--Steinberg Theorem there exists $x \in \rZ^{\circ}(\bL_1)$ such that $(xl)^{-1}F(xl) \in \rZ(\bG)$. If $z = xl$ then $(\bL_1,\delta_1) = ({}^{zg}\bL_2,{}^{zg}\delta_2)$ with $zgZ(\bG) \in \bG_{\ad}^F$ so the map is injective.
\end{proof}

\begin{pa}
After \cref{thm:jor-dec-cusp} we are reduced to understanding the set $\mathcal{J}_{\cusp}^{\circ}(\bG,F)/\bG$. For this we define $\Jor_{\cusp}^{\circ}(\bG,F)$ to be the set of all tuples $((s,n),(\bM,\psi))$ with $(s,n) \in \mathfrak{T}(\bG,F)$ and $(\bM,\psi) \in \Cusp(\rC_{\bG}^{\circ}(s)^{nF})$ a \emph{unipotent cuspidal pair}, meaning that $\psi \in \UCh(\bM,nF)$ is a cuspidal unipotent character. Clearly we have an action of $\bG$ on $\Jor_{\cusp}^{\circ}(\bG,F)$ by setting $g\cdot (\mathfrak{s},(\bM,\psi)) = (g\cdot\mathfrak{s},({}^g\bM,{}^g\psi))$.

The set $\Jor_{\cusp}^{\circ}(\bG,F)/\bG$ parameterises the Harish-Chandra series of unipotent characters and is completely understood via Lusztig's explicit description of all cuspidal unipotent characters, which one finds in \cite[Part 3]{lusztig:1978:representations-of-finite-chevalley-groups} or \cite[\S13.7]{carter:1993:finite-groups-of-lie-type}. The following relates $\Jor_{\cusp}^{\circ}(\bG,F)/\bG$ with $\mathcal{J}_{\cusp}^{\circ}(\bG,F)/\bG$.
\end{pa}

\begin{prop}\label{prop:map-cents-cusp}
We have a well-defined bijection
\begin{align*}
\mathcal{J}_{\cusp}^{\circ}(\bG,F)/\bG &\to \Jor_{\cusp}^{\circ}(\bG,F)/\bG\\
[(s,n),(\bL,\psi)] &\mapsto [(s,n),(\rC_{\bL}^{\circ}(s),\psi)].
\end{align*}
\end{prop}

\begin{proof}
That this is well defined is obvious. Now assume we have pairs $((s_1,1),(\bL_1,\psi_1))$ and $((s_2,1),(\bL_2,\psi_2))$ such that $((s_1,1),(\rC_{\bL_1}^{\circ}(s_1),\psi_1)) = g \cdot ((s_2,1),(\rC_{\bL_2}^{\circ}(s_2),\psi_2))$ for some $g \in \bG$. Then $g \in \bG^F$ and $\rZ_F^{\circ}(\bL_1) = \rZ_F^{\circ}(C_{\bL_1}^{\circ}(s_1)) = {}^g\rZ_F^{\circ}(C_{\bL_2}^{\circ}(s_2)) = {}^g\rZ_F^{\circ}(\bL_2)$. Hence, by \cref{prop:borel-tits} we have $\bL_1 = \rC_{\bG}(\rZ_F^{\circ}(\bL_1)) = {}^g\rC_{\bG}(\rZ_F^{\circ}(\bL_2)) = {}^g\bL_2$. Thus the map is injective.

Suppose $((s,n),(\bM,\psi)) \in \Jor_{\cusp}^{\circ}(\bG,F)$ and let $\bL = \rC_{\bG}(\rZ_{nF}^{\circ}(\bM))$. Note that $\bM \leqslant \bL$ contains a maximal torus of $\rC_{\bG}^{\circ}(s)$ which in turn must contain $s$. By \cref{prop:borel-tits} we have $\bL$ is the smallest $(\bG,F)$-split Levi subgroup containing $\bM = \rC_{\rC_{\bG}^{\circ}(s)}(\rZ_F^{\circ}(\bM)) = \bL \cap \rC_{\bG}^{\circ}(s)$, where here we use that $\bM$ is a $(\rC_{\bG}^{\circ}(s),F)$-split Levi subgroup. As $\bM = \bL \cap \rC_{\bG}^{\circ}(s)$ is connected it follows that $\bM = \rC_{\bL}^{\circ}(s)$ so the map is surjective.
\end{proof}

\begin{rem}
It is straightforward to reduce the classification of the set $\Jor_{\cusp}^{\circ}(\bG,F)/\bG$ to the case where $\bG$ is an adjoint simple group. We refer the reader to \cite[\S11.5]{digne-michel:2020:representations-of-finite-groups-of-lie-type} for more details.
\end{rem}

\begin{pa}\label{pa:extraction-process}
Let us now explain how we can use these results to explicitly determine the set of cuspidal pairs $\Cusp(\bG^F)/\bG_{\ad}^F$ from $\Jor_{\cusp}^{\circ}(\bG^{\star},F^{\star})/\bG^{\star}$. We start with a set of representatives for the orbits $\mathcal{T}_{\bG^{\star}}(\bT_0^{\star},F^{\star})/W^{\star}$ which parameterise the semisimple conjugacy classes of the dual group $\bG^{\star F^{\star}}$ by \cref{lem:ss-reps}. Let $(s,a) \in \mathcal{T}_{\bG^{\star}}(\bT_0^{\star},F^{\star})$. To choose a representative of the coset $a$ we first fix a set of simple roots $\Delta^{\star}(s)$ for the root system $\Phi^{\star}(s) \subseteq \Phi^{\star}$ of $\rC_{\bG^{\star}}^{\circ}(s)$ with respect to $\bT_0^{\star}$. There is then a unique element of the coset $w \in a$ such that $\Delta^{\star}(s)$ is $nF^{\star}$-stable where $n \in \rN_{\bG^{\star}}(\bT_0^{\star})$ is a representative of $w \in W^{\star}$.

Each $nF^{\star}$-stable subset $J \subseteq \Delta(s)$ gives rise to a standard $(\rC_{\bG^{\star}}^{\circ}(s),nF^{\star})$-split Levi subgroup $\bT_0^{\star} \leqslant \bM_J \leqslant \rC_{\bG^{\star}}^{\circ}(s)$ as in \cref{pa:standard-parabolic}. Using Lusztig's classification we list all unipotent cuspidal pairs of the form $(\bM_J,\psi) \in \Cusp(\rC_{\bG^{\star}}^{\circ}(s)^{nF^{\star}})$ with $J \subseteq \Delta(s)$. These give representatives for the $\rC_{\bG^{\star}}^{\circ}(s)^{nF^{\star}}$-orbits of unipotent cuspidal pairs. Moreover, if two such pairs $(\bM_{J_1},\psi_1)$ and $(\bM_{J_2},\psi_2)$ are in the same $\rC_{\bG^{\star}}(s)^{nF^{\star}}$ orbit then $J_1 = J_2$ and $\psi_1 = \psi_2$, see \cite[8.2.1]{lusztig:1984:characters-of-reductive-groups} and the discussion in \cite[\S4]{geck:2018:a-first-guide}. In this way we get representatives $((s,n),(\bM_J,\psi))$ for the orbits $\Jor_{\cusp}^{\circ}(\bG^{\star},F^{\star})/\bG^{\star}$.

Now fix such a tuple and let $\Phi^{\star}(s)_J \subseteq \Phi^{\star}(s)$ be the Levi subsystem determined by $J$ as in \cref{pa:standard-parabolic}. Using \cref{cor:split-levi-cover} we can calculate the smallest $nF^{\star}$-split Levi subsystem $\Sigma^{\star} \subseteq \Phi^{\star}$ containing $\Phi^{\star}(s)_J$. By the proof of \cref{prop:borel-tits} we thus have $\bL^{\star} = \rC_{\bG^{\star}}(\rZ_{nF^{\star}}^{\circ}(\bM_J)) \geqslant \bT_0^{\star}$ is the $(\bG^{\star},nF^{\star})$-split Levi subgroup with root system $\Sigma^{\star}$. The proof of \cref{prop:map-cents-cusp} now shows that $\bM_J = \rC_{\bL^{\star}}^{\circ}(s)$.

We have thus obtained a tuple $((s,n),(\bL^{\star},\psi)) \in \mathcal{J}_{\cusp}^{\circ}(\bG^{\star},F^{\star})$ and this gives a set of representatives for the $\bG^{\star}$-orbits. It follows from \cref{lem:conj-to-standard} that there exists an $\rN_{\bG^{\star}}(\bT_0^{\star})$-conjugate $((s',n'),(\bL_I^{\star},\psi'))$ of this tuple with $\bL_I^{\star}$ a $(\bG^{\star},F^{\star})$-split standard Levi subgroup which is dual to the $(\bG,F)$-split standard Levi subgroup $\bL_I$ with $I \subseteq \Delta$. If $\delta \in \Irr(\bL_I^F)$ is such that $J^{\bL_I}(\delta) = [(s',n'),\psi']$ then $(\bL_I,\delta) \in $ is a cuspidal pair and we have inverted the map in \cref{thm:jor-dec-cusp}.
\end{pa}

\begin{rem}
Assume $F$ is a Frobenius endomorphism. Then the above can be easily generalised to the $d$-cuspidal setting as follows. For an integer $d \geqslant 1$ we define $\Cusp_d(\bG^F)$, $\mathcal{J}_{d,\cusp}^{\circ}(\bG,F)$, and $\Jor_{d,\cusp}^{\circ}(\bG,F)$, in the same way except replacing ``cuspidal'' by ``$d$-cuspidal'' and ``$(\bG,F)$-split'' by ``$d$-split'', see \cite[\S3.5]{geck-malle:2020:the-character-theory-of-finite-groups-of-lie-type}. The map in \cref{thm:jor-dec-cusp} then gives an injection $\Cusp_d(\bG^F)/\bG_{\ad}^F \to \mathcal{J}_{d,\cusp}^{\circ}(\bG,F)/\bG$ and the map in \cref{prop:map-cents-cusp} gives a bijection $\mathcal{J}_{d,\cusp}^{\circ}(\bG,F)/\bG \to \Jor_{d,\cusp}^{\circ}(\bG,F)/\bG$. Again, in this case we can invert the map in \cref{prop:map-cents-cusp} using the results of \cref{sec:split-parabolics}.
\end{rem}

\section{Cuspidal Pairs in Isolated Lusztig Series}\label{sec:iso-series}
\begin{pa}\label{pa:simplex}
In this section we assume that $\bG$ is an adjoint simple group and that $F$ is a Frobenius endomorphism. This assumption means that the automorphism $\phi \in \GL(X_{\mathbb{R}})$ defined in \cref{pa:reductive-setup} stabilises the lattice $X$ and also the sets of roots $\Delta \subseteq \Phi^+ \subseteq \Phi$. We will follow the procedure described in \cref{pa:extraction-process} to obtain representatives for the orbits $\mathcal{J}_{\cusp}^{\circ}(\bG,F)/\bG$.

Let $\mathbb{Z}_{(p)}$ be the localisation of $\mathbb{Z}$ at the prime ideal $(p)$. If $\iota : \mathbb{Z}_{(p)}/\mathbb{Z} \to \mathbb{K}^{\times}$ is a fixed isomorphism then we denote by $\tilde{\iota} : \mathbb{Q} \to \mathbb{K}^{\times}$ the homomorphism obtained as the composition
\begin{equation*}
\mathbb{Q} \to \mathbb{Q}/\mathbb{Z} \to \mathbb{Z}_{(p)}/\mathbb{Z} \overset{\iota}{\to} \mathbb{K}^{\times}.
\end{equation*}
Here the first map is the natural projection map and the second map is obtained by quotienting out the $p$-torsion subgroup of $\mathbb{Q/Z}$.

For brevity we let $V = X_{\mathbb{Q}}$ and dually we let $\wc{V} = \mathbb{Q}\otimes_{\mathbb{Z}} \wc{X}$. We then have a unique surjective homomorphism of abelian groups $\tilde{\iota}_0 : \wc{V} \to \bT_0$ satisfying $\tilde{\iota}_0(k\otimes\gamma) = \gamma(\tilde{\iota}(k))$. We denote by $\wc{\Omega} = (\wc{\varpi}_{\alpha})_{\alpha \in \Delta}$ the basis of $\wc{V}$ dual to the simple roots $\Delta \subseteq V$ with respect to $\langle-,-\rangle : X\times \wc{X} \to \mathbb{Z}$ (extended linearly). As we assume $\bG$ is adjoint we have $X = \mathbb{Z}\Delta$ and $\wc{X} = \mathbb{Z}\wc{\Omega}$.

Recall that $s \in \bT_0$ is said to be \emph{$\bG$-isolated} if $\rC_{\bG}^{\circ}(s)$ is not contained in any proper Levi subgroup of $\bG$. If $\Phi(s) = \{\alpha \in \Phi \mid \alpha(s) = 1\}$ then $s$ is $\bG$-isolated if and only if $\Phi(s)$ is $\Phi$-isolated, see \cref{def:Phi-isolated}. Denote by $\alpha_0 = \sum_{\alpha \in \Delta} n_{\alpha}\alpha \in \Phi$ the highest root and let $\wt{\Delta} = \Delta\cup\{-\alpha_0\}$ be the extended set of simple roots. By convention we set $\wc{\varpi}_{-\alpha_0} = 0$ and $n_{-\alpha_0} = 1$. If $\wt{\Delta}_{p'}$ is the set of $\alpha \in \wt{\Delta}$ such that $\gcd(n_{\alpha},p) = 1$ then the following is shown in \cite[Thm.~5.1]{bonnafe:2005:quasi-isolated}.
\end{pa}

\begin{prop}
Every $\bG$-isolated element in $\bT_0$ is $W$-conjugate to some $h_{\alpha} := \tilde{\iota}_0(\wc{\varpi}_{\alpha}/n_{\alpha})$ with $\alpha \in \wt{\Delta}_{p'}$. Moreover, the root system $\Phi(h_{\alpha})$ of the centraliser $\rC_{\bG}^{\circ}(h_{\alpha})$ with respect to $\bT_0$ has a simple system given by $\Delta(h_{\alpha}) := \wt{\Delta}\setminus\{\alpha\}$.
\end{prop}

\begin{pa}
To understand the orbits $[h_{\alpha},a] \in \mathcal{T}_{\bG}(\bT_0,F)/W$ it suffices to understand the orbits of $\rC_W(h_{\alpha})$ acting on the set of cosets $\mathcal{A}_W(h_{\alpha},F)$. However, after \cite[Thm.~5.1]{bonnafe:2005:quasi-isolated} this is easy once we find an element $w \in W$ satisfying ${}^wF(h_{\alpha}) = h_{\alpha}$. As we will now explain this is straightforward except for a few cases in $\E_8$. If $\alpha \in \wt{\Delta}_{p'}$ let $w_{\alpha} \in \rC_W^{\circ}(h_{\alpha})$ be the unique element satisfying ${}^{w_{\alpha}}\Delta(h_{\alpha}) = -\Delta(h_{\alpha})$ so that $w_0 := w_{-\alpha_0}$ is the longest element of $W$ with respect to $\Delta$.
\end{pa}

\begin{lem}
If $\alpha \in \wt{\Delta}_{p'}$ then either $F(h_{\alpha}) = h_{\alpha}$ or ${}^{w_0}F(h_{\alpha}) = h_{\alpha}$ unless $n_{\alpha} = 5$, in which case $\bG$ is of type $\E_8$, and $q \equiv 2,3 \pmod{5}$.
\end{lem}

\begin{proof}
A straightforward calculation shows that $F(h_{\alpha}) = h_{\phi^{-1}\alpha}^q$ and ${}^{w_0}h_{\alpha} = h_{\varepsilon^{-1}\alpha}^{-1}$ where $\varepsilon = -w_0 \in \GL(X)$ satisfies $\varepsilon\Delta = \Delta$. From the classification of irreducible root systems we see that $1 \leqslant n_{\alpha} \leqslant 6$. Note $n_{\alpha}$ is the order of $h_{\alpha}$ so if $n_{\alpha} = 1$ then certainly $F(h_{\alpha}) = h_{\alpha}$. Assume $n_{\alpha} > 2$ then one easily checks that $\alpha$ is fixed under all automorphisms of the Dynkin diagram. Hence, if $q \equiv 1 \pmod{n_{\alpha}}$ then $F(h_{\alpha}) = h_{\alpha}$ and if $q \equiv -1 \pmod{n_{\alpha}}$ then ${}^{w_0}F(h_{\alpha}) = h_{\alpha}$. As $\gcd(q,n_{\alpha}) = 1$ one of these conditions must hold unless $n_{\alpha} = 5$. Finally suppose $n_{\alpha} = 2$. Then either $\alpha$ is fixed by $\phi$ or $\bG$ is of type $\E_6$ and $\phi = \varepsilon$. In either case the statement holds as $\varepsilon$ has order $2$.
\end{proof}

\begin{pa}
We will now determine the tuples $[(s,n),(\bL,\psi)] \in \mathcal{J}_{\cusp}^{\circ}(\bG,F)/\bG$ when $s$ is $\bG$-isolated and $\bG$ is of classical type, see \cref{pa:exceptional-example} for groups of exceptional type. For these purposes we consider the case of ${}^3\D_4$ to be of exceptional type. The relevant information is given in \cref{tab:HC-series-clsc-typs}. We omit the unipotent character $\psi$ from the table because, according to Lusztig's classification \cite[Part 3]{lusztig:1978:representations-of-finite-chevalley-groups}, in all cases the group $\rC_{\bL}^{\circ}(s)^{nF}$ admits a unique cuspidal unipotent character.

As in \cref{pa:extraction-process} we start with a pair $(h_{\alpha},a) \in \mathcal{T}_{\bG}(\bT_0,F)$, with $\alpha \in \tilde{\Delta}_{p'}$, and choose an element $n \in \rN_{\bG}(\bT_0)$ such that $\overline{n} = w \in a$ is the unique element of the coset satisfying ${}^{w\phi}\Delta(h_{\alpha}) = \Delta(h_{\alpha})$.
\end{pa}

\begin{rem}
If ${}^{w_0}F(h_{\alpha}) = h_{\alpha}$ then $w = w_{\alpha}w_0$ is the unique element of the coset $\rC_W^{\circ}(h_{\alpha})w_0$ satisfying ${}^{w\phi}\Delta(h_{\alpha}) = \Delta(h_{\alpha})$.
\end{rem}

\begin{pa}
We now list all the unipotent cuspidal pairs $(\bM_J,\psi) \in \Cusp(\rC_{\bG}^{\circ}(h_{\alpha}),nF)$ with $J \subseteq \Delta(h_{\alpha})$. This gives the column entitled $\rC_{\bL}^{\circ}(s)^{nF}$ of \cref{tab:HC-series-clsc-typs}. We let $\theta : V \to V^{w\phi}$ be the projection onto the $w\phi$-fixed space of $V$. For brevity we denote by $V_1 \subseteq V$ the subspace spanned by the subsystem $\Phi(h_{\alpha})_J$. Our aim is to calculate the smallest $w\phi$-split Levi subsystem $\Sigma = \Phi \cap \theta^{-1}(\theta(V_1))$ containing $\Phi(h_{\alpha})_J$. We will specify $\Sigma$ by specifying the set of positive roots $\Sigma^+ = \Sigma \cap \Phi^+$.
\end{pa}

\subsection{Classical Types}
\begin{pa}
The case of type $\A_n$ is trivial, as the only isolated element is the identity, so we assume that $\bG$ has type $\B_n$ ($n \geqslant 2$), $\C_n$ ($n \geqslant 2$), or $\D_n$ ($n \geqslant 4$). We take the root system and simple roots $\Delta = (\alpha_1,\dots,\alpha_n)$ to be labelled as in \cite{bourbaki:2002:lie-groups-chap-4-6}. We will assume that $(e_1,\dots,e_n)$ is a basis of the $\mathbb{Q}$-vector space $V$ so that the roots are expressed as in \cite{bourbaki:2002:lie-groups-chap-4-6} with respect to this basis.

Let $\mathcal{E} = \{e_1,\dots,e_n,-e_n,\dots,-e_1\} \subseteq V$. The natural action of $W$ on $V$ preserves the set $\mathcal{E}$ so we have a permutation representation $W \to \mathfrak{S}_{\mathcal{E}}$ into the symmetric group on $\mathcal{E}$. This representation is faithful and we will identify $W$ with its image. For $1 \leqslant i < j \leqslant n$ we define the following elements of $\mathfrak{S}_{\mathcal{E}}$: $p_{i,j} = p_{j,i} = (i,j)(-i,-j)$ and $u_{i,j} = u_{j,i} = (i,-j)(-i,j)$. Moreover, for $1 \leqslant i \leqslant n$ we define $c_i = (i,-i)$.
\end{pa}

\begin{pa}
We now consider a few of the more interesting cases from \cref{tab:HC-series-clsc-typs}. We identify the case by listing the type of $\rC_{\bG}^{\circ}(s)^{nF}$. In all the special cases we consider we have $\phi$ is the identity. Before proceeding we make a few comments on the conventions used in \cref{tab:HC-series-clsc-typs}. Here we use the following notation for the sets of square and triangular integers:
\begin{itemize}
	\item $\square_m = \{(2k+m)^2 \mid k \in \mathbb{Z}_{\geqslant 0}\}\setminus\{1\}$ so that $\square_0 = \{0,4,16,\dots\}$ and $\square_1 = \{9,25,49,\dots\}$,
	\item $\triangle = \{k(k+1)/2 \mid k \in \mathbb{Z}_{\geqslant 0}\} = \{0,1,3,6,10,15,\dots\}$.
\end{itemize}
These integers occur in the classification of cuspidal unipotent characters for groups of classical type, see \cite[Thm.~8.2]{lusztig:1977:irreducible-representations-of-finite-classical-groups}.

Suppose the root system $\Phi(s)$ of $\rC_{\bG}^{\circ}(s)$ has two indecomposable components $\mathsf{X}_e\cdot \mathsf{Y}_{n-e}$. In this case a component of $\rC_{\bL}^{\circ}(s)$ with rank $d \leqslant e$ is taken to be contained in the component $\mathsf{X}_e$ and a component with rank $m-d \leqslant n-e$ or $m \leqslant n$ is taken to be contained in the component $\mathsf{Y}_{n-e}$. Finally, let us note that each row of the table corresponds to exactly one orbit of pairs $(s,a) \in \mathcal{T}_{\bG}(\bT_0,F)$ except the rows marked by ($\star$) which correspond to two such orbits.
\end{pa}

\begin{pa}[Case $\C_{\frac{n}{2}}(q^2)$]\label{pa:case-Cn2}
We take $w = u_{1,n}u_{2,n-1}\cdots u_{n/2,n/2+1}$. The fixed point space $V^w$ has as basis $(v_i \mid 1 \leqslant i \leqslant \frac{n}{2})$ with $v_i = e_i-e_{n+1-i}$. Moreover, the subspace $\theta(V_1) \subseteq V^w$ has as basis $(v_1,\dots,v_{\frac{m}{2}})$. If $K = \{1,\dots,\frac{m}{2},n+1-\frac{m}{2},\dots,n\}$ then $\Sigma^+$ consists of $\{2e_i,e_i \pm e_j \mid i,j \in K$ and $i < j\}$, giving a component of type $\C_m$, and $\{e_i+e_{n+1-i} \mid \frac{m}{2} < i \leqslant \frac{n}{2}\}$ giving the $\frac{n-m}{2}$ components of type $\A_1$.
\end{pa}

\begin{pa}[Case ${}^2\D_e(q)\cdot\B_{n-e}(q)$]
Recall that $2 \leqslant e \leqslant n$. We take $w = c_1$. The fixed point space $V^w$ has as basis $(e_2,e_3,\dots,e_n)$. There are two cases to consider. If $\bM_J$ has type $\B_m$ then $\theta(V_{\Phi_1})$ has as basis $(e_{n-m+1},\dots,e_n)$. If $K = \{1,n-m+1,\dots,n\}$ then $\Sigma^+$ is given by $\{e_i,e_i\pm e_j \mid i,j \in K$ and $i < j\}$ which is of type $\B_{m+1}$.

If $\bM_J$ has type $\D_d\cdot\B_{m-d}$ then $\theta(V_{\Phi_1})$ has as basis $(e_1,\dots,e_d,e_{n-m+d+1},\dots,e_n)$. If $K = \{1,\dots,d,n-m+d+1,\dots,n\}$ then $\Sigma^+$ is given by $\{e_i,e_i\pm e_j \mid i,j \in K$ and $i < j\}$ which is of type $\B_m$.
\end{pa}

\begin{pa}[Case ${}^2\D_e(q)\cdot{}^2\D_{n-e}(q)$]
Recall that $2 \leqslant e \leqslant \frac{n}{2}$. We take $w = c_1c_n$. The fixed point space $V^w$ has as basis $(e_2,e_3,\dots,e_{n-1})$. First, if $\bM_J = \bT_0$ is a torus then $\Sigma^+ = \{e_1-e_n,e_1+e_n\}$, which has type $\D_2 = \A_1\A_1$.

If $\bM_J$ has type $\D_m$ then $\theta(V_{\Phi_1})$ has as basis $(e_{n-m+1},\dots,e_{n-1})$. If $K = \{1,n-m+1,n-m+2,\dots,n\}$ then $\Sigma^+$ is given by $\{e_i\pm e_j \mid i,j \in K$ and $i < j\}$ which is of type $\D_{m+1}$. The case where $\bM_J$ is of type $\D_d$ or $\D_d\cdot \D_{m-d}$ is similar to previous cases. 
\end{pa}

\begin{pa}[Case $\D_{\frac{n}{2}}(q^2)$]
Here there are two orbits of pairs $(s,a)$ such that $\rC_{\bG}^{\circ}(s)^{nF}$ has type $\D_{\frac{n}{2}}(q^2)$. We take $w$ to be one of the following two elements:
\begin{equation*}
\begin{cases}
u_{1,n}u_{2,n-1}\cdots u_{n/2,n/2+1}\\
p_{1,n}u_{2,n-1}\cdots u_{n/2,n/2+1}
\end{cases}
\end{equation*}
We have $V^w$ has as basis $(v_i \mid 1 \leqslant i \leqslant \frac{n}{2})$ with $v_i = e_i-e_{n+1-i}$ for $i > 1$ and $v_1 = e_1-e_n$ for the first choice of $w$ and $v_1 = e_1+e_n$ for the second choice of $w$. We now proceed as in \cref{pa:case-Cn2}.
\end{pa}

\subsection*{The Exceptional Types}
\begin{pa}\label{pa:exceptional-example}
If $\bG$ is of exceptional type then the analogous data to that found in \cref{tab:HC-series-clsc-typs}, amongst other things, can be extracted from the tables of \cite{kessar-malle:2013:quasi-isolated-blocks}. This data can also be easily recomputed using \cite{taylor:2020:cusppairs} and commands provided in \textsf{CHEVIE}. Specifically, in \cite{taylor:2020:cusppairs} we have implemented the command \texttt{JordanCuspidalLevis} which can be used to obtain the relevant Levi subgroups. We give an example in the case of $\E_8$, see \cite[Table 7, no.~19--24]{kessar-malle:2013:quasi-isolated-blocks}
\end{pa}

\begin{verbatim}
gap> GF := CoxeterCoset(CoxeterGroup("E",8), ());;
gap> s := SemisimpleElement(Group(GF), [0,0,0,0,0,1/3,0,0]);;
gap> CG := Centralizer(Group(GF), s).group;;
gap> w := LongestCoxeterElement(CG)*LongestCoxeterElement(G);;
gap> CGF := CoxeterSubCoset(GF, InclusionGens(CG), w);
2E6<1,5,3,4,2,96>x2A2
gap> JordanCuspidalLevis(CGF);
[ [ E8, 2E6<1,5,3,4,2,96>x2A2 ],
  [ E7<6,2,5,4,3,1,74>.(q-1), 2E6<1,5,3,4,2,96>.(q-1)(q+1) ],
  [ E7<1,2,3,4,13,7,8>.(q-1), 2A5<1,3,4,2,96>x2A2<7,8>.(q-1) ],
  [ D6<2,13,4,3,1,74>.(q-1)^2, 2A5<1,3,4,2,96>.(q-1)^2(q+1) ],
  [ D4<26,8,7,27>.(q-1)^4, 2A2<7,8>.(q-1)^4(q+1)^2 ],
  [ A1<26>xA1<27>xA1<74>.(q-1)^5, (q-1)^5(q+1)^3 ] ]
\end{verbatim}

\begin{landscape}
\begingroup
\small
\begin{table}
\caption{Isolated Elements in Classical Types.}\label{tab:HC-series-clsc-typs}
\centering
\begin{tabular}{>{$}c<{$}*{5}{>{$}l<{$}}}
\toprule
\multicolumn{1}{c}{$\bG^F$} & \multicolumn{1}{c}{$C_{\bG}^{\circ}(s)^{nF}$} & p? & \multicolumn{1}{c}{$\bL^{nF}$} & \multicolumn{1}{c}{$C_{\bL}^{\circ}(s)^{nF}$} & \multicolumn{1}{c}{Conditions}\\
\midrule

\B_n(q) & \B_n(q) & - & \B_m(q)\cdot\Phi_1^{n-m} & \B_m(q)\cdot\Phi_1^{n-m} & m \in 2\triangle\\
& \D_e(q)\cdot\B_{n-e}(q)\quad (2 \leqslant e \leqslant n) & \neq 2 & \B_m(q)\cdot\Phi_1^{n-m} & \D_d(q)\cdot\B_{m-d}(q)\cdot\Phi_1^{n-m} & d \in \square_0, m-d \in 2\triangle\\
& {}^2\D_e(q)\cdot\B_{n-e}(q)\quad (2 \leqslant e \leqslant n) & \neq 2 & \B_{m+1}(q)\cdot\Phi_1^{n-m-1} & \B_m(q)\cdot\Phi_1^{n-m-1}\cdot\Phi_2 & m \in 2\triangle\\
& & & \B_m(q)\cdot\Phi_1^{n-m} & {}^2\D_d(q)\cdot\B_{m-d}(q)\cdot\Phi_1^{n-m} & d \in \square_1, m-d \in 2\triangle\\
\midrule

\C_n(q) & \C_n(q) & - & \C_m(q)\cdot\Phi_1^{n-m} & \C_m(q)\cdot\Phi_1^{n-m} & m \in 2\triangle \\
& \C_e(q)\cdot\C_{n-e}(q)\quad (1 \leqslant e \leqslant n/2) & \neq 2 & \C_m(q)\cdot\Phi_1^{n-m} & \C_d(q)\cdot\C_{m-d}(q)\cdot\Phi_1^{n-m} & d, m-d \in 2\triangle\\
& \C_{\frac{n}{2}}(q^2) \quad (n\text{ even}) & \neq 2 & \A_1(q)^{\frac{n-m}{2}}\cdot\C_m(q)\cdot\Phi_1^{\frac{n-m}{2}} & \C_{\frac{m}{2}}(q^2)\cdot\Phi_1^{\frac{n-m}{2}}\cdot\Phi_2^{\frac{n-m}{2}} & m \in 4\triangle\\
\midrule

\D_n(q) & \D_n(q) & - & \D_m(q)\cdot\Phi_1^{n-m} & \D_m(q)\cdot\Phi_1^{n-m} & m \in \square_0\\
& \D_e(q)\cdot\D_{n-e}(q)\quad (2 \leqslant e \leqslant n/2) & \neq 2 & \D_m(q)\cdot\Phi_1^{n-m} & \D_d(q)\cdot\D_{m-d}(q)\cdot\Phi_1^{n-m} & d, m-d \in \square_0\\
& {}^2\D_e(q)\cdot{}^2\D_{n-e}(q)\quad (2 \leqslant e \leqslant n/2) & \neq 2 & \D_2(q)\cdot\Phi_1^{n-2} & \Phi_1^{n-2}\cdot\Phi_2^2 & \\
& & & \D_{d+1}(q)\cdot\Phi_1^{n-d-1} & {}^2\D_d(q)\cdot\Phi_1^{n-d-1}\cdot\Phi_2 & d \in \square_1\\
& & & \D_{m+1}(q)\cdot\Phi_1^{n-m-1} & {}^2\D_m(q)\cdot\Phi_1^{n-m-1}\cdot\Phi_2 & m \in \square_1\\
& & & \D_m(q)\cdot\Phi_1^{n-m} & {}^2\D_d(q)\cdot{}^2\D_{m-d}(q)\cdot\Phi_1^{n-m} & d,m-d \in \square_1\\
(\star) & \D_{\frac{n}{2}}(q^2)\quad (n\text{ even}) & \neq 2 & \A_1(q)^{\frac{n-m}{2}}\cdot\D_m(q)\cdot\Phi_1^{\frac{n-m}{2}} & \D_{\frac{m}{2}}(q^2)\cdot\Phi_1^{\frac{n-m}{2}}\cdot\Phi_2^{\frac{n-m}{2}} & m \in 2\square_0\\
\midrule

{}^2\D_n(q) & {}^2\D_n(q) & - & \Phi_1^{n-1}\cdot\Phi_2 & \Phi_1^{n-1}\cdot\Phi_2 & \\
& & & {}^2\D_m(q)\cdot\Phi_1^{n-m} & {}^2\D_m(q)\cdot\Phi_1^{n-m} & m \in \square_1\\
& \D_e(q)\cdot{}^2\D_{n-e}(q)\quad (2 \leqslant e \leqslant n-2) & \neq 2 & \Phi_1^{n-1}\cdot\Phi_2 & \Phi_1^{n-1}\cdot\Phi_2 & \\
& & & {}^2\D_m(q)\cdot\Phi_1^{n-m} & \D_d(q)\cdot{}^2\D_{m-d}(q)\cdot\Phi_1^{n-m} & d \in \square_0, m-d \in \square_1 \\
& {}^2\D_e(q)\cdot\D_{n-e}(q)\quad (2 \leqslant e \leqslant n-2) & \neq 2 & \Phi_1^{n-1}\cdot\Phi_2 & \Phi_1^{n-1}\cdot\Phi_2 & \\
& & & {}^2\D_{m+1}(q)\cdot\Phi_1^{n-m-1} & \D_m(q)\cdot\Phi_1^{n-m-1}\cdot\Phi_2 & m \in \square_2\\
& & & {}^2\D_m(q)\cdot\Phi_1^{n-m} & {}^2\D_d(q)\cdot \D_{m-d}(q)\cdot\Phi_1^{n-m} & d \in \square_1, m-d \in \square_0\\
(\star) & {}^2\D_{\frac{n}{2}}(q^2)\quad (n\text{ even}) & \neq 2 & \A_1(q)^{\frac{n-2}{2}}\cdot{}^2\D_2(q)\cdot\Phi_1^{\frac{n-2}{2}} & \Phi_1^{\frac{n-2}{2}}\cdot\Phi_2^{\frac{n-2}{2}}\cdot\Phi_4 & \\
& & & \A_1(q)^{\frac{n-m}{2}}\cdot{}^2\D_m(q)\cdot\Phi_1^{\frac{n-m}{2}} & {}^2\D_{\frac{m}{2}}(q^2)\cdot\Phi_1^{\frac{n-m}{2}} & m \in 2\square_1\\
\bottomrule
\end{tabular}
\end{table}
\endgroup
\end{landscape}

\setstretch{0.96}
\renewcommand*{\bibfont}{\small}
\printbibliography
\end{document}